\documentclass[12pt]{amsart}  %% {article}
\usepackage{amsmath,amsfonts,amsthm}

\def\Tau{{\mathcal T}}

\newcommand{\R}{\mathbb{R}}
\newcommand{\grad}{\hbox{\rm grad}}

\newcommand{\Half}{{\textstyle{\frac{1}{2}}}}
\newcommand{\BigOh}[1]{\mathcal{O}(#1)}

\numberwithin{equation}{section}
\numberwithin{theorem}{section}

\begin{document}

\title[On the Hamiltonian for water waves]
      {On the Hamiltonian for water waves}

\author{Walter Craig}
%% \author{Catherine Sulem}
\address{Department of Mathematics, McMaster University, Hamilton ON,
  L8S 4K1 Canada, \ craig@math.mcmaster.ca}
%% \address{Department of Mathematics, University of Toronto, Toronto ON,
%%   M5S 2E4 Canada,  \ sulem@math.toronto.edu}
\thanks{This article is an extension of a lecture at the RIMS
  Symposium on Mathematical Analysis in Fluid and Gas Dynamics, July 6
  - 8 2016. The author's work is partially supported by the Canada Research
  Chairs Program and NSERC through grant number 238452--16, and he is grateful for the
  hospitality of the Courant Institute during the preparation of this article.}
%% \thanks{C.S. is  partially  supported by NSERC through grant
%%  number 46179--13}


\date{\today}

\begin{abstract}
Many equations that arise in a physical context can be posed in the form 
of a Hamiltonian system, meaning that there is a symplectic structure on 
an appropriate phase space, and a Hamiltonian functional with respect to
which time evolution of their solutions can be expressed in terms of a
Hamiltonian vector field.  It is known from the work of VE Zakharov
that the equations for water waves can be posed as a Hamiltonian
dynamical system, and that the equilibrium solution is an elliptic
stationary point. In this article we generalize the Hamiltonian
formulation of water waves by Zakharov to a general coordinatization
of the dynamical free surface, which allows it to apply to situations
that include overturning wave profiles. This answers a question posed
to the author by T.~Nishida during the RIMS Symposium on Mathematical
Analysis in Fluid and Gas Dynamics that took place during July 6 - 8 2016.  
\end{abstract}

\maketitle

\section{Introduction}

The equations for water waves describe the flow of an incompressible
and irrotational fluid with a free surface, under the additional
restoring forces of gravity and with the possibility to include the
effects of surface tension.  The fluid velocity $u(t,x,y)$, expressed 
in Eulerian coordinates, satisfies the conditions
\begin{equation}\label{Eqn:IncompressibleIrrotational}
   \nabla \cdot u = 0 ~, \qquad \nabla \wedge u = 0 ~,
\end{equation}
in a fluid domain $\Omega(t) \subseteq \R^{d-1}_x \times \R^1_y$ whose
boundary consists of two components, a bottom described by a
hypersurface $s \in \R^{d-1} \mapsto b(s) \in \R^d$ and a free
surface given by a time dependent hypersurface 
$s \in \R^{d-1} \mapsto \gamma(t,s) \in \R^d$. %% \{ (x,y) = \gamma(t,s)\}$. 
It is possible that the bottom is unbounded below, and indeed it is
common to consider the case that the bottom boundary lies at 
$\{ y = -\infty \}$. 
The free surface of the fluid domain itself, defined by the
hypersurface $\gamma(t,s)$, is one of the unknowns. Because of the
constraints \eqref{Eqn:IncompressibleIrrotational} the fluid motion is
given by a {\em potential flow}  
\looseness=-10000
\begin{align}\label{Eqn:Laplace}
   & u = \nabla \varphi ~, \qquad \Delta \varphi = 0 ~, \quad 
     \hbox{\rm in the fluid domain} \ \Omega(t) ~;  \\
   & \partial_N \varphi = 0 ~, \qquad \hbox{\rm bottom boundary
       conditions on} \ (x,y) \in \{ b(s) \ : s \in \R^{d-1}\} ~.   \nonumber
\end{align}
Denote the horizontal and vertical components of the hypersurface
defining the free surface of $\Omega(t)$ by 
$\gamma(t,s) = (\gamma_1(t,s), \gamma_2(t,s))$ (that is, 
$x = \gamma_1 \in \R^{d-1}$ and $y = \gamma_2 \in \R^1$), and
the space-time unit normal vector to the free surface to be
${\bf N}_{t,x,y}$. Furthermore define the space-time vector describing the
fluid velocity field as 
${\bf T}_{t,x,y} = (1,u(t,x,y))^T = (1,\nabla_{x,y}\varphi)^T$. 
Then the kinematic free boundary condition on the free surface 
$\{(x,y) \in \gamma\}$ is the geometrical condition that
\begin{equation}\label{Eqn:KinematicBC}
    {\bf N}_{t,x,y} \cdot {\bf T}_{t,x,y} = 0 ~.
\end{equation}
The physics of the flow of Euler's equations is described by the
second nonlinear boundary condition;  
\begin{eqnarray}\label{Eqn:FreeSurfaceBC}
     \partial_t\varphi & = & -g\gamma_2 - \Half |\nabla \varphi |^2 
   + \sigma H(\eta)  \qquad
      \ \hbox{\rm Bernoulli condition} ~. 
\end{eqnarray}
The force of surface tension is given by the term $\sigma H$, where
$H(\eta)$ is the mean curvature of the free surface. We study both of
the cases $\sigma > 0$ and $\sigma = 0$. 

In the case that the the free surface is given as a graph, $y = \eta(t,x)$,
$x \in \R^{d-1}$ the conditions \eqref{Eqn:KinematicBC}\eqref{Eqn:FreeSurfaceBC} 
can be rewritten as
\begin{eqnarray}\label{Eqn:GraphFreeSurfaceBC}
    \partial_t \eta & = & \partial_y \varphi - \partial_x \eta
     \cdot \partial_x \varphi    \qquad\qquad \hbox{\rm kinematic boundary conditions}   \\
     \partial_t\varphi & = & -g\eta -\Half |\nabla \varphi |^2 
   + \sigma H(\eta)  \qquad
      \ \hbox{\rm Bernoulli condition} ~. \nonumber
\end{eqnarray}
This dynamic free boundary problem was recognized  by VE~Zakharov
(\cite{Z68}~1968) to be a Hamiltonian PDE, which is to say that equations 
\eqref{Eqn:Laplace}\eqref{Eqn:GraphFreeSurfaceBC} can be given the form of
a Hamiltonian system 
\[
    \dot z = X^H(z) ~, \qquad \hbox{\rm where} \quad X^H(z) = J\, \grad_z  H(z) ~,
\]
with the Hamiltonian function $H$ the total energy of the system
\eqref{Eqn:IncompressibleIrrotational}\eqref{Eqn:Laplace}\eqref{Eqn:GraphFreeSurfaceBC}, 
namely 
\begin{equation}\label{Eqn:Energy}
   H = \frac{1}{2} \iint_{\Omega(t)} |\nabla \varphi|^2 \, dydx 
     + {\frac{g}{2}}\int_{\R^{d-1}} \eta^2
   \, dx + \sigma \int_{\R^{d-1}} \sqrt{1 + |\partial_x \eta|^2} - 1 \, dx ~. 
\end{equation}
A more subtle aspect is the choice of canonical variables for the phase 
space, which as in ~\cite{Z68} is normally given by 
$z(x) := (\eta(x), \xi(x) := \varphi(x,\eta(x)))$, namely
\begin{align*}
   & \partial_t\begin{pmatrix}\eta \\ \xi \end{pmatrix} 
     = \begin{pmatrix} 0 & I \\ -I & 0 \end{pmatrix} 
       \begin{pmatrix} \delta_\eta H \\ \delta_\xi H \end{pmatrix} ~, 
     \qquad \begin{pmatrix} 0 & I \\ -I & 0 \end{pmatrix} := J ~,
\end{align*}
which is in Darboux coordinates in that the symplectic form is defined by
$J$ as given above.

The goals of this article are to explain this fact, and to extend
the formulation of the equations of water waves as a Hamiltonian
system to free surfaces given in general coordinates, satisfying equations
\eqref{Eqn:IncompressibleIrrotational}\eqref{Eqn:Laplace} and boundary
conditions \eqref{Eqn:KinematicBC}\eqref{Eqn:FreeSurfaceBC}. In
particular this allows the equations to describe the evolution of
smooth overturning wave profiles.  

\section{Hamiltonian for overturning waves}

During the RIMS Symposium on Mathematical analysis in fluid and gas
dynamics in Kyoto, T. Nishida asked the author whether Zakharov's
formulation of the water waves problem \eqref{Eqn:GraphFreeSurfaceBC}
as a Hamiltonian PDE could be extended to take into consideration the
case of geometries of free surfaces that are not graphs, and in particular
waves that are overturning. A subtle issue in Zakharov's formulation is the
specific choice of canonical conjugate variables, which appears to
require remarkable insight, but in retrospect can be deduced from a 
{\em principle of least action} \`a la Lagrange and a subsequent
Legendre transform \cite{C08}. It turns out that similar
considerations are useful when seeking to describe the water wave
problem in general coordinates.  

%%%***

After a first version of this article was circulated, Tom Bridges
brought to my attention the article of T.B.~Benjamin and P.~Olver~\cite{BO82},
which gives a `quasi-Hamiltonian' structure for the problem of water
waves with arbitrary parametrization of the free surface. It also is a
formulation that can describe overturning wave profiles. The
result was elaborated in T.~Bridges \& N.~Donaldson~\cite{BD11}. The
principal difference between the formulation in these two articles and
the present one, in addition to being derived for both two- and
three-dimensional cases, is that in the former the symplectic form is
degenerate due to the extra degrees of freedom of the parametrization
of the free surface, and it depends upon the phase space variables, in
particular on gradient of the velocity potential. In the present article
Hamilton's canonical equations are nondegenerate, given in Darboux
coordinates, and by using the Dirichlet -- Neumann operator the
evolution equations can be restricted entirely to the free surface. 

We will address Nishida's question essentially on a formal level, and in
the case $d=2$, for which $\Omega(t) \subseteq \R^2$. Configuration
space is taken to be the space of curves 
$\Gamma := \{s \mapsto \gamma(s) \ : s \in \R^1 \}$. To actually perform
analysis for data in this configuration space, including solving
Laplace's equation on the fluid domain $\Omega(t)$, we
should give some topology to this space, such as $\gamma \in C^1(\R^1)$, 
and we should consider free surfaces $\gamma(s)$ that have a limit
$\lim_{s\to \pm\infty} \gamma(s) = 0$. Furthermore we should ask that
there be a uniform lower bound on the distance between $\gamma$ and
the bottom boundary $\{b(s)\}$, and also that $\gamma$ satisfy
a global chord - arc condition. However in the present context we will
ignore these details.  

\subsection{Free surface boundary conditions}

Given a one parameter family of curves $\gamma(t,s)$, the velocity, which is
the time derivative $\partial_t\gamma(t,s) = \dot\gamma(t,s)$ defines
a vector field in the tangent space over the curve $\gamma$. A
natural orthonormal frame for the tangent space over $\gamma(s)$ is
given by $(T(s),N(s))$, where   
\[
   T(s) = \frac{\partial_s\gamma(s)}{|\partial_s\gamma(s)|} ~, \quad 
   N(s) = -JT(s) ~, \quad 
   J = \begin{pmatrix} 0 & 1 \\ -1 & 0 \end{pmatrix} ~.
\]
In the frame $(T(s),N(s))$ the velocity $\dot\gamma$ vector field
can be represented by its coordinates
\[
   n(t,s) = N\cdot\dot\gamma ~, \quad \tau(t,s) = T\cdot\dot\gamma ~.
\]
The equations for water waves determine the evolution of the fluid
domain $\Omega(t)$ and the velocity field $u(t,x,y)$ defined in 
$\Omega(t)$. Because of the condition of irrotationality 
\eqref{Eqn:IncompressibleIrrotational} they can be reduced to two
nonlinear boundary conditions posed on the free surface
$\gamma(t,s)$. The first of these is the kinematic condition 
\eqref{Eqn:KinematicBC}, which in this frame is written
\[
    0 = {\bf N}_{t,x,y} \cdot {\bf T}_{t,x,y} 
      = c(t,s) [N\cdot \dot\gamma - N\cdot u]~, 
\]
where $c(t,s)$ is a normalization irrelevant to the present
discussion. When interpreted in terms of the velocity potential
$\nabla_{x,y}\varphi(t,x,y) = u(t,x,y)$ this is the statement that
\begin{equation}\label{Eqn:Kinematic-2}
   N\cdot\nabla\varphi\bigl|_\gamma = N\cdot\dot\gamma = n(t,s) ~.
\end{equation} 
With the foresight of Zakharov's formulation, define 
$\xi(s) = \varphi(\gamma(s))$ to be the boundary values of the velocity
potential on the free surface $\gamma(s)$. Then the function $n(s)$
can be expressed in terms of the Dirichlet -- Neumann operator for the
fluid domain
\begin{equation}\label{Eqn:n-intermsof-p}
     n(s) = G(\gamma)\xi(s) ~
\end{equation}
(and these quantites will depend parametrically on time $t$).
Namely, $\varphi$ is the solution of Laplace's equation on the fluid
domain $\Omega$ satisfying Neumann boundary conditions on the bottom 
$(x,y) = b(s)$ and with boundary data on the free surface 
$\gamma(s)$ given by $\varphi(\gamma(s)) = \xi(s)$, where the 
operator $G(\gamma)$ is then defined by
\[
   \xi(s) \mapsto \varphi(x,y) \mapsto \nabla\varphi\cdot N :=
   G(\gamma)\xi(s) ~. 
\] 
The normalization for the Dirichlet -- Neumann
operator is that $|N| = 1$ (which differs slightly from what is
commonly used for the problem posed in graph coordinates). This is an
elliptic boundary value problem which can be solved for
$\varphi(x,y)$, hence the map  $\xi \mapsto G(\gamma)\xi$ is well defined. 
 
The Bernoulli condition \eqref{Eqn:FreeSurfaceBC} expresses the
physics described by the Euler equations on the free surface. Written
in terms of $\xi(t,s) = \varphi(\gamma(t,s))$, for which 
$\partial_t\xi(t,s) = \partial_t\varphi(t,\gamma(t,s)) 
   = \varphi_t + \nabla\varphi \cdot \dot\gamma$, 
this is 
\begin{equation}\label{Eqn:BernoulliCond-p}
   \partial_t \xi - \nabla\varphi\cdot\gamma = -g\gamma_2 - \Half|\nabla\varphi|^2~.
\end{equation}
Recalling the definition that $\dot\gamma = nN(t,s) + \tau T(t,s)$
(and using that $|T| = |N| = 1$ and $T\cdot N = 0$),
\[
     \nabla\varphi\cdot\dot\gamma = (\nabla\varphi \cdot N)n +
     (\nabla\varphi \cdot T) \tau = (G(\gamma)\xi) n 
     + \frac{\partial_s \xi}{|\partial_s\gamma|} \tau ~. 
\]
Therefore \eqref{Eqn:BernoulliCond-p} is rewritten as
\begin{equation}\label{Eqn:BernoulliCond-p2}
   \partial_t \xi = -g\gamma_2 
    + \Half\bigl[(G(\gamma)\xi)^2 -
      \frac{1}{|\partial_s\gamma|^2}(\partial_s \xi)^2 
    + 2\frac{\partial_s \xi}{|\partial_s\gamma|} \tau \bigr] ~,
\end{equation}
where we have used the definition \eqref{Eqn:n-intermsof-p} for $n$ in
terms of $\xi$.  

In general the geometry of the curve $\gamma$ can be recovered from
$T(s)$ (or equivalently from $N(s)$) by integration, but not its 
parametrization. However so far in this discussion we have not
addressed the issues of ambiguity that have been introduced by
allowing arbitrary (nonsingular) coordinatization of curves
$\gamma(s)$. There exist numerous useful possibilities to specify this
parametrization, one standard one being as a graph, but another is to
parametrize by arc length. In this latter case  
\begin{align*}
    & \partial_s\gamma = T ~, \quad |T(s)| = 1 ~, \\
    & \partial_s T = \kappa(s) N ~, \quad \partial_s N = -\kappa(s) T
  ~,
\end{align*}
which describes the evolution in $s$ of the Frenet frame, $\kappa(s)$
being the curvature. In these coordinates one recovers $\tau(s)$ from 
$n(s)$; indeed because
\[
    0 = \partial_t |T(t,s)|^2 = 2 \partial_t T  \cdot T
      = 2 \partial_s\dot\gamma \cdot T 
\]
then one has
\begin{equation}\label{Eqn:TangentalVelocity}
   \partial_s\tau = \partial_s\dot\gamma\cdot T +
   \dot\gamma\cdot\partial_s T = \kappa \dot\gamma\cdot N = \kappa n ~.  
\end{equation}
In arc length coordinates, equation \eqref{Eqn:BernoulliCond-p2} is
somewhat simpler, namely
\begin{equation}\label{Eqn:BernoulliCond-p3}
   \partial_t \xi = -g\gamma_2 + \Half\bigl[(G(\gamma)\xi)^2 
     - (\partial_s \xi)^2 + 2\partial_s \xi \tau \bigr] ~.
\end{equation}
In this case, the tangential component of the velocity is recovered
from \eqref{Eqn:n-intermsof-p}\eqref{Eqn:TangentalVelocity}, namely 
$\partial_s\tau = \kappa n = \kappa G(\gamma)\xi$.

\subsection{Legendre transform}

The Lagrangian for free surface water waves corresponds to the total
energy of the system, which consists of two terms, the kinetic energy
$K$ and the potential energy $U$;
\begin{equation}\label{Eqn:Lagrangian}
   L = K - U ~.
\end{equation}
The Legendre transform is the classical approach to transfer a
Lagrangian system into the canonical conjugate coordinates of a
Hamiltonian system. When the Lagrangian functional $L$ is expressed in
terms of the variables $\gamma$ and $\dot\gamma$, by analogy with
classical mechanics one defines {\em conjugate momentum variables} via
the Legendre transform as
\[
   \xi = \delta_{\dot\gamma} L ~.
\]

The kinetic energy is given by the Dirichlet integral
\[
   K = \iint_{\Omega} \Half|\nabla \varphi(x,y)|^2 \, dydx
\] 
and the potential energy is respectively
\[
     U = \iint_{\Omega} gy \, dydx + C ~,
\]
which is, as usual, only defined up to an additive
constant. If the effects of surface tension were to be included in the
equations of motion, then the potential energy has an additional term,
namely
\[
    U = \iint_{\Omega} gy \, dydx + \sigma\int_\gamma \, dS_\gamma + C' ~,
\]
where $dS_\gamma = |\partial_s\gamma(s)| \, ds$. 
Our derivation below is in the case that $\sigma = 0$, but by
modifications of the argument the case $\sigma \not= 0$ is also able
to be included. 

Integrating by parts in $K$ and using the boundary
conditions, we can express the kinetic energy in terms of 
integrated quantities on the free surface
\begin{equation}\label{Eqn:KineticEnergy}
    K = \int_{\gamma} \Half \xi G(\gamma)\xi \, dS_{\gamma} ~.
\end{equation}
We note that the normalization for the Dirichlet -- Neumann operator
$G(\gamma)$ is different from that used in \cite{Z68} and \cite{CS93},
so that it is Hermetian with respect to the line element $dS_{\gamma}$. 
Using \eqref{Eqn:n-intermsof-p} the kinetic energy can be written in
terms of $\gamma$ and $\dot\gamma$;
\begin{equation}\label{Eqn:KineticEnergy-2}
     K(\gamma,\dot\gamma) := \int_{\gamma} \Half n G^{-1}(\gamma) n \,
     dS_{\gamma} ~. 
\end{equation}
The potential energy $U$ can be expressed with respect to the divergence
theorem, using a vector field $V(x,y) := (0,\frac{g}{2}y^2)^T$;
\begin{align}\label{Eqn:PotentialEnergy-2}
   U(\gamma) & = \iint_{\Omega} \nabla \cdot V(x,y) \, dvol 
       = \int_{\gamma} V \cdot N \, dS_{\gamma} + C 
    = \int_{\gamma} \frac{g}{2} \gamma_2^2 
       \frac{\partial_s\gamma_1}{|\partial_s\gamma|} \, dS_{\gamma} 
       + C ~. 
\end{align}
In arc length parametrization this would read
\[
    U(\gamma) = \int_{\gamma} \frac{g}{2} \gamma_2^2 
       \partial_s\gamma_1 \, ds + C ~. 
\]

In the case of general coordinates for the free surface
$\gamma(s)= (\gamma_1(s),\gamma_2(s))$, gradients are expressed with
respect to the metric given by $\int_\gamma \cdot \, dS_\gamma$. 
The tangent space $T_\gamma$ at $\gamma$ to the set of curves is given
coordinates using the Frenet frame $(T(s),N(s))$. Variations of $L$ with
respect to vector fields $Y(s) \in T_\gamma$ along $\gamma(s)$ can be
decomposed into their normal and tangential components, namely 
\[
    \langle \delta_{Y} L , \delta Y\rangle_\gamma 
     = \int_{\gamma} \grad_{N \cdot Y} L \, (N \cdot \delta Y)    
     + \grad_{T \cdot Y} L \, (T \cdot \delta Y) ~dS_{\gamma}~.
\]
If the vector field $Y(s) = \dot\gamma(t,s)$ is the velocity of the
curve $\gamma(t,s)$, this is written
\[
    \langle \delta_{\dot\gamma} L, \delta\dot\gamma \rangle_\gamma  
     = \int_{\gamma} \grad_{N \cdot \dot\gamma} L \, (N \cdot \delta\dot\gamma)    
     + \grad_{T \cdot \dot\gamma} L \, (T \cdot \delta\dot\gamma) ~dS_{\gamma}~.
\]
In the case of the kinetic energy $K$ above, and because of the
decomposition $\dot\gamma(s) = \tau(s) T(s) + n(s) N(s)$, this is
\begin{equation}\label{Eqn:LegendreTransform}
    \delta_{\dot\gamma} K = \delta_n K + \delta_\tau K 
     = G^{-1}(\gamma)n + 0 ~.
\end{equation}
Thus $\xi = G^{-1}(\gamma)n$ is the canonical conjugate variable to
normal perturbations of a given free surface $\gamma$, while $\tau$
remains undefined without further specification of the parametrization
of the curve $\gamma$. This degeneracy will be resolved when a
particular form of parametrization is imposed. 

Following the prescription of the Legendre transform 
\eqref{Eqn:LegendreTransform} the Hamiltonian is given by
\begin{equation}
   H = K + U = \frac{1}{2} \int_{\gamma} \xi G(\gamma)\xi \, dS_{\gamma} 
     + \frac{g}{2} \int_{\gamma}  
       \gamma_2^2\frac{\partial_s\gamma_1}{|\partial_s\gamma|} \, dS_{\gamma} ~.
\end{equation}
The remaining questions are to how to best express the variables that
are canonically conjugate to $\xi(s)$, and to show that the resulting
equations of motion \eqref{Eqn:n-intermsof-p}\eqref{Eqn:BernoulliCond-p2}  
coincide with the Hamiltonian vector field, namely
\begin{align*}\label{Eqn:HamiltonianSystem} 
   \partial_t z = J \grad H(z) ~. 
\end{align*}

The gradient of the kinetic energy $K$ with respect to $\xi$ is 
\[
   \grad_\xi K = G(\gamma) \xi ~, 
\] 
which corresponds to the conjugate of the normal variations of $K$ with
respect to $\dot\gamma$. 

Using the expression \eqref{Eqn:PotentialEnergy-2}, 
the gradient of the potential energy is given by
\begin{align}\label{Eqn:GradientPotentialEnergy}
   \langle \delta_\gamma U, \delta\gamma\rangle_\gamma & = \int_\gamma \frac{g}{2} 
     \begin{pmatrix} -2\gamma_2\partial_s\gamma_2 \\
                      2\gamma_2\partial_s\gamma_1 \end{pmatrix}
     \cdot \begin{pmatrix} \delta\gamma_1
       \\ \delta\gamma_2 \end{pmatrix}
     \, ds \\
   & = \int g \gamma_2(s) \, N \cdot \begin{pmatrix} \delta\gamma_1
       \\ \delta\gamma_2 \end{pmatrix} \, dS_\gamma  ~,  \nonumber 
\end{align} 
corresponding to the gradient of $U$ with respect to normal variation
of $\gamma$ itself, namely $\grad_{N\cdot\delta\gamma} U$. 

The gradient of the kinetic energy $K$ with respect to $\gamma$ is the
more subtle quantity in this formulation. Consider a fluid domain
$\Omega$ with free surface $\gamma(s)$ and a family of nearby domains 
$\Omega_1$ with nearby free surfaces 
$\gamma_1(s) = \gamma(s) + \delta\gamma(s)$. Denote the outward unit
normal by $N(s)$ and $N_1(s)$ respectively. We consider the Dirichlet
integrals
\[
   K(\gamma,\xi) = \Half\int_\gamma \xi(s) G(\gamma)\xi(s) \, dS_\gamma ~,
   \qquad
   K_1 = K(\gamma_1,\xi) = \Half\int_{\gamma_1} \xi(s) G(\gamma_1)\xi(s) \, dS_{\gamma_1} ~,
\]
for which we impose that the boundary values of the velocity potentials
$\Phi_1(x,y)$ on $\gamma_1$ and $\Phi(x,y)$ on $\gamma$ coincide
\[
    \Phi(\gamma(s)) = \xi(s) = \Phi_1(\gamma_1(s)) ~,
\]
while we vary the boundary curve $\gamma(s)$ to 
$\gamma_1(s) = \gamma(s) + \delta\gamma(s)$. This is to say that one
takes the partial derivative of the kinetic energy with respect to
variations of the domain, while fixing the boundary conditions for the
velocity potential on the free surface. To this effect, the
boundary values of $\Phi(x,y)$ on the curve $\gamma_1(s)$ are given by
\begin{align*}
    \Phi (\gamma_1(s)) & = \Phi(\gamma(s)) + \nabla\Phi(\gamma(s)) \cdot
    \delta\gamma(s) + \BigOh{\delta^2} \\
    & = \Phi(\gamma(s)) + (\nabla\Phi\cdot N) \, N\cdot \delta\gamma(s) 
      + (\nabla\Phi\cdot T) \, T \cdot \delta\gamma(s) + \BigOh{\delta^2}
      ~. 
\end{align*}
Therefore 
\begin{equation}\label{Eqn:VariationPhi}
    \Phi_1(\gamma_1(s)) - \Phi(\gamma_1(s)) = - (\nabla\Phi\cdot N) \, N\cdot \delta\gamma(s) 
      - (\nabla\Phi\cdot T) \, T \cdot \delta\gamma(s) + \BigOh{\delta^2}
\end{equation}
Furthermore, given a harmonic function $\Phi(x,y)$ defined on a
neighborhood that includes $\Omega \cup \Omega_1$, by Green's theorem
the difference of the boundary integral expressions for their
Dirichlet integrals is given by  
\begin{align}\label{Eqn:VariationPhi-2}
 & \Half \int_{\gamma_1} \Phi(\gamma_1(s)) N_1\cdot \nabla \Phi(\gamma_1(s))
   \, dS_{\gamma_1} 
   - \Half \int_{\gamma} \Phi(\gamma(s)) N\cdot \nabla \Phi(\gamma(s))
   \, dS_{\gamma} \\ 
  & = \Half \iint_{\Omega_1 \backslash \Omega} 
    |\nabla \Phi|^2 \, dvol 
   \simeq \Half \int_\gamma |\nabla \Phi|^2 N\cdot \delta\gamma(s) \, dS_\gamma ~.  
     \nonumber
\end{align}
Therefore the variation of the kinetic energy $K$ with fixed boundary
data $\xi(s)$ is calculated as the limit in small $\delta$ of 
\begin{align*}
   K_1 - K = & \Half\int_{\gamma_1} \xi(s) G(\gamma_1)\xi(s) \, dS_{\gamma_1} 
      - \Half\int_{\gamma} \xi(s) G(\gamma)\xi(s) \, dS_{\gamma} \\ 
     = & \Half\int_{\gamma_1} \Phi_1(\gamma_1) N_1\cdot
      \nabla\Phi_1(\gamma_1) \, dS_{\gamma_1} 
     - \Half\int_{\gamma} \Phi(\gamma) N\cdot
      \nabla\Phi(\gamma) \, dS_{\gamma} \\ 
     = & \int_{\gamma_1} (\Phi_1 - \Phi)(\gamma_1) N_1\cdot
      \nabla\Phi_1(\gamma_1) \, dS_{\gamma_1}   \\
     & + \Half\int_{\gamma_1} \Phi(\gamma_1)
     N_1\cdot\nabla\Phi(\gamma_1) \, dS_{\gamma_1} 
     - \Half\int_{\gamma} \Phi(\gamma)
     N\cdot\nabla\Phi(\gamma) \, dS_{\gamma} + \BigOh{\delta^2} ~. 
\end{align*}
Using \eqref{Eqn:VariationPhi} in the first term and \eqref{Eqn:VariationPhi-2} 
in the second and third, 
\begin{align*}
   K_1 - K = & \int_\gamma - (\nabla\Phi\cdot N)^2 \, N\cdot \delta\gamma(s) 
      - (\nabla\Phi\cdot N)(\nabla\Phi\cdot T) \, T \cdot \delta\gamma(s) 
      dS_\gamma \\
      & + \Half \int_\gamma |\nabla\Phi|^2 N\cdot \delta\gamma \,
      dS_\gamma + \BigOh{\delta^2} ~. 
\end{align*}
Furthermore, both of the velocity potentials $\Phi$ and $\Phi_1$
satisfy Neumann boundary conditions on the bottom $(x,y) = b(s)$.
Thus $N\cdot\nabla\Phi(\gamma(s)) = G(\gamma)\xi(s)$ and 
$T\cdot\nabla\Phi(\gamma(s)) = \frac{1}{|\partial_s\gamma|}\partial_s\xi(s)$, 
giving an expression in the limit as $\delta \to 0$ for $\grad_{\delta\gamma} K$, 
namely  
\begin{align}\label{Eqn:GradK}
   &  \langle \delta K \cdot \delta\gamma \rangle_\gamma 
     = \int_\gamma \grad_{\delta\gamma} K \cdot \delta\gamma \, dS_\gamma 
   \\
     & \quad = \Half \int_\gamma - \bigl(G(\gamma)\xi \bigr)^2 N\cdot \delta\gamma 
     + \Bigl(\frac{1}{|\partial_s\gamma|} \partial_s \xi \Bigr)^2 N\cdot \delta\gamma
     - 2\Bigl(\frac{1}{|\partial_s\gamma|} \partial_s \xi G(\gamma)\xi \Bigr) 
       T\cdot \delta\gamma \, dS_\gamma ~.   \nonumber
\end{align}
With these expressions in hand, we conclude that the equations of
motion for the problem of water waves takes the canonical form of a
Hamiltonian system; 
\begin{align}
  & N\cdot \partial_t\gamma = \grad_\xi H \\
  & \partial_t \xi = - \grad_{N\cdot\delta\gamma} H ~. \nonumber
\end{align}

In general the choice of coordinatization of the free surface is made
separately from the decomposition of the tangent space $T_\gamma$ into
its normal and tangential components. Variations $Y(s) = \delta\gamma(s)$ 
of $\gamma$ are necessarily constrained by the coordinate choice to
the class of {\em admissible variations}. The choice of
coordinitization determines the tangential component of the velocity
$\tau = T\cdot\partial_t\gamma$ as a function of the normal component,
through the constraints imposed by the coordinatization of the free surface. 
This applies in particular to the time derivative of the curve,
$\dot\gamma(s) \in T_\gamma$. That is, coordinitization dictates a relation
between $T\cdot\delta\gamma$ and $N\cdot \delta\gamma$, say
$T\cdot\delta\gamma = \Tau(\gamma) (N\cdot\delta\gamma)$ 
in somewhat abstract terms. Thus, in terms of such a coordinate choice,
\begin{equation}\label{Eqn:GradientK}
   \grad_{N\cdot\delta\gamma} K = 
     \Half\Bigl[\Bigl(\frac{1}{|\partial_s\gamma|} \partial_s \xi \Bigr)^2 
     - \bigl(G(\gamma)\xi \bigr)^2  
     - 2\Bigl(\frac{1}{|\partial_s\gamma|} 
        \partial_s \xi G(\gamma)\xi \Tau(\gamma)\Bigr) \Bigr] ~.
\end{equation}
This gradient is worked out in detail for several
standard choices of parametrization in the subsection below. 

\subsection{Particular coordinates}

Common choices for the parametrization of the free surface are:
(1) the classical case of free surfaces given as a graph in 
$x \in \R^1$, which does not allow for overturning free surfaces. 
(2) arc length parametrization of $\gamma(s)$ which are able
to describe overturning wave profiles. In these coordinates we have
seen that $\partial_s \tau = \kappa n$. 
(3) Lagrangian coordinates, for which fluid particle positions are
advected by the flow, $\partial_t(X(t),Y(t)) = u(X(t),Y(t)) 
   = \nabla\varphi(\gamma(t,\cdot))$, or 
(4) conformal mapping coordinates as used in \cite{KN79}. Specifying
the coordinatization of free surface curves $\gamma$ in cases (1)(2) and
(4) gives rise to systems of constraints which may be considered to be
holonomic as they are imposed independently of the velocity $\dot\gamma$. 
The parametric specification by Lagrangian coordinates in contrast is a
nonholonomic constraint. 

The traditional choice of parametrization is  (1) to write the surface as a
graph; in such graph coordinates, where $\gamma = (x,\eta(x))$, 
and the pair of variables $(\eta(x),\xi(x))$ are canonically conjugate as
given by Zakharov~\cite{Z68}. With the expression for the kinetic
energy $K$ in terms of the Dirichlet -- Neumann operator $G(\eta)$ as
in~\cite{CS93}, then  
\[
    H(\eta,\xi) = \int_{\R^1} \Half \xi G(\eta) \xi \, \sqrt{1 +
      (\partial_x\eta)^2} \, dx  
     + \frac{g}{2} \int_{\R^1} \eta^2 \, dx ~.
\]
In these graph coordinates, $\gamma = (x,\eta(x))$ so that admissible
variations are $\delta\gamma = (0,\delta\eta)$, and the relationship
between $N\cdot\delta\gamma$ and $T\cdot\delta\gamma$ is given by
\[
    T\cdot\delta\gamma = \partial_x \eta \, N\cdot\delta\gamma ~.
\]
The gradient of the kinetic energy is thus
\[
  \Half\Bigl[ \Bigl(\frac{(\partial_x \xi)^2}{1 + (\partial_x\eta)^2} \Bigr) 
     - \bigl(G(\eta)\xi \bigr)^2 
     - 2\Bigl(\frac{1}{\sqrt{1 + (\partial_x\eta)^2}} 
       \partial_x \xi G(\eta)\xi \, \partial_x\eta\Bigr) \Bigr] ~.
\]
Because of this, 
\[
  N\cdot\dot\gamma = \frac{1}{\sqrt{1 + (\partial_x\eta)^2}} 
   \begin{pmatrix} -\partial_x \eta \\ 1 \end{pmatrix} 
     \cdot \begin{pmatrix} 0 \\ \partial_t \eta \end{pmatrix} 
\]
and the resulting equations \eqref{Eqn:GraphFreeSurfaceBC} for water
waves are given by 
\begin{align*}  
   & \frac{1}{\sqrt{1 + (\partial_x\eta)^2}}\partial_t \eta 
      = \delta_\xi H = G(\eta)\xi \\ 
   & \partial_t \xi = -\delta_\gamma H = -g\eta 
      + \Half\Bigl[(G(\eta)\xi)^2 
      - \frac{(\partial_x \xi)^2}{1 + (\partial_x\eta)^2}
      + \frac{2\partial_x \xi  G(\eta)\xi \, \partial_x\eta}{\sqrt{1 + (\partial_x\eta)^2}}
        \Bigr] ~. 
\end{align*}
Calculating for an independent verification of \eqref{Eqn:BernoulliCond-p2}, 
one finds that in graph coordinates 
\[
   \tau = \frac{\partial_t\eta\partial_x\eta}
               {\sqrt{1 + (\partial_x\eta)^2}} 
        = G(\eta)\xi \, \partial_x\eta ~. 
\]
This system of equations, modulo the difference in normalization of
the Dirichlet -- Neumann operator $G(\eta)$, appears in \cite{CS93},
and is used in the existence theory for water waves and many of its
distinguished scaling limits in~\cite{L05}\cite{L13}.

Coordinates given in terms of arc length along the free surface
$\gamma(s)$ allow the system \eqref{Eqn:KinematicBC}\eqref{Eqn:FreeSurfaceBC} 
to describe overturning wave profiles. This choice of coordinates
implies in particular that $\partial_s\gamma(s) = T(s)$ and
$\partial_s T \perp T$, since
\[
    |\partial_s\gamma(s)|^2 = 1 ~,  \quad 
    0 = \partial_s |\partial_s\gamma(s)|^2 
      = 2\partial_s\gamma \cdot \partial_s^2\gamma ~. 
\]
Indeed any vector field $Y(s)$ along the curve $\gamma(s)$ that arises
from an infinitessimal motion which preserves the arc length
parametrization must satisfy
\[
    0 = \frac{d}{d\delta}\Bigl|_{\delta=0}|\partial_s\gamma + \delta Y|^2 =
    2\partial_s\gamma\cdot Y = 2T\cdot Y ~. 
\]
Admissible variations $\delta\gamma(s)$ are arc length preserving in
the present case, implying that $Y = \partial_s\delta\gamma$ is as
above, and hence  
\[
     0 = \partial_s(T\cdot\delta\gamma) = \partial_s T\cdot\delta\gamma +
     T\cdot \partial_s\delta\gamma = \kappa N\cdot\delta\gamma ~. 
\]
This is the relationship between tangential and normal variations that
applies to the gradient of the kinetic energy, an interesting
geometrical aspect of this choice of coordinates. The resulting
Bernoulli equations of motion are 
\[
   \partial_t \xi = -g\gamma_2 - \Half\Bigl[
     \Bigl(\frac{1}{|\partial_s\gamma|} \partial_s \xi \Bigr)^2 
     - \bigl(G(\gamma)\xi \bigr)^2  
     - 2\Bigl(\frac{1}{|\partial_s\gamma|} 
       \partial_s \xi G(\gamma)\xi \Tau(\gamma)\Bigr) \Bigr]
\]
where $\Tau(\gamma)$ satisfies 
\[
    \partial_s\bigl(G(\gamma)\xi \Tau\bigr) = \kappa G(\gamma)\xi ~.
\]

%%%%%%%%%%%%%%%%%%%%%%%%%%%%%%%%%%%%%%%%%%%%%%%%%%%%%%%%%%%%%%%%%%%
%%%%%%%%%%%%%%%%%%%%%%%%%%%%%%%%%%%%%%%%%%%%%%%%%%%%%%%%%%%%%%%%%%%

\end{document}